\documentclass[12pt]{article}
\usepackage{amsmath,amssymb,amsthm}
\usepackage{fullpage}

\newtheorem{thm}{Theorem}

\newcommand{\dis}{\displaystyle}

\begin{document}

\title{Ramanujan's Harmonic Number Expansion}
\author{Mark B. Villarino\\
Depto.\ de Matem\'atica, Universidad de Costa Rica,\\
2060 San Jos\'e, Costa Rica}
\date{\today}

\maketitle

 \begin{abstract}
 An algebraic transformation of the \textsc{DeTemple-Wang} half-integer approximation to the harmonic series produces the general formula and error estimate for the \textsc{Ramanujan} expansion for the $n$th harmonic number.\end{abstract}

\section{Introduction}
Entry~9 of Chapter~38 of B. \textsc{Berndt}'s edition of
\textsc{Ramanujan}'s Notebooks, Volume 5 \cite[p.~521]{Berndt} reads:
\begin{quote}
``\textit{Let $m :=  \frac{n(n+1)}{2}$, where $n$ is a positive
integer. Then, as $n$ approaches infinity,}
\begin{align*}
\sum_{k=1}^n \frac{1}{k}
\sim \frac{1}{2} \ln(2m) + \gamma + \frac{1}{12m} - \frac{1}{120m^2}
+ \frac{1}{630m^3} - \frac{1}{1680m^4} + \frac{1}{2310m^5} \\
-\frac{191}{360360m^6} + \frac{29}{30030m^7} - \frac{2833}{1166880m^8}+\frac{140051}{17459442m^{9}} -[\cdots]
\text{.''}\\
\end{align*}
\end{quote}

Berndt's proof simply verifies (as he himself explicitly notes) that
Ramanujan's expansion coincides with the standard \textsc{Euler}
expansion
\begin{align*}
H_n := \sum_{k=1}^n \frac{1}{k}
&\sim \ln n + \gamma + \frac{1}{2n} - \frac{1}{12n^2}
+ \frac{1}{120n^4} -[ \cdots]
\\
&= \ln n + \gamma - \sum_{k=1}^\infty \frac{B_k}{n^k} 
\end{align*}
where $B_k$ denotes the $k^{\mathrm{th}}$ \textsc{Bernoulli} number
and $\gamma := 0.57721\cdots$ is Euler's constant.

However, Berndt does \textit{not} give the \emph{general formula} for the coefficient of $\frac{1}{m^k}$ in  \textsc{Ramanujan}'s expansion, nor does he prove that it is an
\textit{asymptotic} series in the sense that the error in the value obtained
by stopping at any particular stage in \textsc{Ramanujan}'s series is less than
the next term in the series. Indeed we have been unable to find
\emph{any} error analysis of \textsc{Ramanujan}'s series.

Although it is true that the asymptotics of the harmonic numbers were already determined by the \textsc{Euler} expansion, later mathematicians have offered alternative approximative formulas (see \cite{Cesaro, CQ, D, Lodge, TM, Vill2}).  \textsc{Ramanujan}'s formula is one of the most accurate (see \cite{Vill2})..

However, the \textsc{Euler} expansion apparently does not easily lend itself to an error analysis of \textsc{Ramanujan}'s expansion, nor to a general formula for the coefficient of $\frac{1}{m^{k}}.$

In an earlier paper (see \cite{Vill}) we proved that the first five terms of the \textsc{Ramanujan} expansion indeed form an asymptotic expansion in the sense above, but we did not prove the general case.

The general case is the subject of the following:

\begin{thm}
For any integer $p\geqslant 1$ define:
\begin{equation}
\fbox{$\dis R_{p}:=\frac{(-1)^{p-1}}{2p\cdot8^{p}}\left\{1+\sum_{k=1}^{p}\binom{p}{k}(-4)^{k}B_{2k}\left(\frac{1}{2}\right)\right\}$},
\end{equation} where $B_{2k}(x)$ is the \textsc{Bernoulli} polynomial of order $2k$.  Put\begin{align}
\fbox{$\dis
m:
= \frac{n(n+1)}{2}
$}
\end{align}where $n$ is a positive integer.
Then, for every integer $r\geqslant 1$, there exists a $\Theta_r$, $0 < \Theta_r < 1$, for which the
following equation is true:
\begin{align}
\fbox{$\dis
1+ \frac{1}{2}+\frac{1}{3}+\cdots+\frac{1}{n}
= \frac{1}{2} \ln(2m) + \gamma + \sum_{p=1}^{r}\frac{R_{p}}{m^{p}}+\Theta_r\cdot\frac{R_{r+1}}{m^{r+1}}
$}
\end{align}
\end{thm}

We observe that the formula for $R_{p}$ can be written symbolically as follows: 
\begin{align}
\fbox{$\dis
R_p
= -\frac{1}{2p}\left(\frac{4B^{2}-1}{8}\right)^{p}
$}
\end{align}where we write $B_{2m}\left(\frac{1}{2}\right)$ in place of $B^{2m}$ after carrying out the above expansion.
%---------------------------------------------------------------------------------------------
\section{Proof of Ramanujan's Expansion}
\begin{proof}

We begin with the half-integer approximation to $H_{n}$ due to \textsc{DeTemple} and \textsc{Wang} (see \cite{D}): \emph{For any positive integer $r$ there exists a $\theta_r$, $0<\theta_r<1,$ for which the following equation is true}:
\begin{align}
H_n=\ln \left(n+\frac{1}{2}\right)+ \gamma + \sum_{p=1}^{r}\frac{D_{p}}{\left(n+\frac{1}{2}\right)^{2p}}+\theta_r\cdot\frac{D_{r+1}}{\left(n+\frac{1}{2}\right)^{2r+2}}    
\end{align}\emph{where} \begin{align}
 D_p &:=-\frac{B_{2p}\left(\frac{1}{2}\right)}{2p}.   
\end{align}Since$$\left(n+\frac{1}{2}\right)^2=2m+\frac{1}{4}$$

\noindent we obtain:\begin{align*}
\sum_{p=1}^{r}\frac{D_{p}}{\left(n+\frac{1}{2}\right)^{2p}}&=\sum_{p=1}^{r}\frac{D_{p}}{(2m)^{p}\left(1+\frac{1}{8m}\right)^{p}}\\
&=\sum_{p=1}^{r}\frac{D_{p}}{(2m)^{p}}\left(1+\frac{1}{8m}\right)^{-p}\\
&=\sum_{p=1}^{r}\frac{D_{p}}{(2m)^{p}}\sum_{k=0}^\infty\binom{-p}{k}\frac{1}{8^k m^k}\\
&=\sum_{p=1}^{r}\frac{D_{p}}{2^{p}}\sum_{k=0}^\infty(-1)^k\binom{k+p-1}{k}\frac{1}{8^k}\cdot\frac{1} {m^{p+k}}\\
&=\sum_{p=1}^{r}\left\{\sum_{s=0}^{p-1}\frac{D_{s}}{2^{s}}(-1)^{p-s}\binom{p-1}{p-s}\frac{1}{8^{p-s}}\right\}\cdot\frac{1} {m^{p}}+E_r\\
\end{align*}where\begin{align*}
      E_r:=&\frac{D_{1}}{2^{1}}\sum_{k=r}^\infty(-1)^k\binom{k}{k}\frac{1}{8^k}\cdot\frac{1} {m^{1+k}}
+\frac{D_{2}}{2^{2}}\sum_{k=r-1}^\infty(-1)^k\binom{k+1}{k}\frac{1}{8^k}\cdot\frac{1} {m^{2+k}} 
+\cdots\\
&+\frac{D_{r}}{2^{r}}\sum_{k=1}^\infty(-1)^k\binom{k+r-1}{k}\frac{1}{8^k}\cdot\frac{1} {m^{r+k}}  
\\
\end{align*}

Substituting the right hand side of the last equation into the right hand side of (5) we obtain:

\begin{align}H_n=\ln \left(n+\frac{1}{2}\right)+ \gamma + \sum_{p=1}^{r}\left\{\sum_{s=0}^{p-1}\frac{D_{s}}{2^{s}}(-1)^{p-s}\binom{p-1}{p-s}\frac{1}{8^{p-s}}\right\}\cdot\frac{1} {m^{p}}+E_r+\theta_r\cdot\frac{D_{r+1}}{\left(n+\frac{1}{2}\right)^{2r+2}}.    
\end{align}

Moreover,
\begin{align*}
\ln \left(n+\frac{1}{2}\right)&=\frac{\ln \left(n+\frac{1}{2}\right)^2}{2}\\
&=\frac{1}{2}\ln\left(2m+\frac{1}{4}\right)\\
&=\frac{1}{2}\ln(2m)+\frac{1}{2}\ln\left(1+\frac{1}{8m}\right)\\
&=\frac{1}{2}\ln(2m)+\frac{1}{2}\sum_{l=1}^{\infty}(-1)^{l-1}\frac{1}{l8^lm^l}\\
&=\frac{1}{2}\ln(2m)+\frac{1}{2}\sum_{l=1}^{r}(-1)^{l-1}\frac{1}{l8^lm^l}+\epsilon_r\\
\end{align*}where $$\epsilon_r:=\sum_{l=r+1}^{\infty}(-1)^{l-1}\frac{1}{2l8^lm^l}.$$

Substituting the right-hand side of this last equation into (7) we obtain\begin{align*}H_n=&\frac{1}{2}\ln(2m)+\frac{1}{2}\sum_{l=1}^{r}(-1)^{l-1}\frac{1}{l8^lm^l}+ \gamma + \sum_{p=1}^{r}\left\{\sum_{s=0}^{p-1}\frac{D_{s}}{2^{s}}(-1)^{p-s}\binom{p-1}{p-s}\frac{1}{8^{p-s}}\right\}\cdot\frac{1} {m^{p}}\\
&+\epsilon_r+E_r+\theta_r\cdot\frac{D_{r+1}}{\left(n+\frac{1}{2}\right)^{2r+2}}]\\
=&\frac{1}{2}\ln(2m)+\gamma+\sum_{p=1}^{r}\left\{(-1)^{p-1}\frac{1}{2p8^p} + \sum_{s=0}^{p-1}\frac{D_{s}}{2^{s}}(-1)^{p-s}\binom{p-1}{p-s}\frac{1}{8^{p-s}}\right\}\cdot\frac{1} {m^{p}}\\
&+\epsilon_r+E_r+\theta_r\cdot\frac{D_{r+1}}{\left(n+\frac{1}{2}\right)^{2r+2}}\\    
\end{align*}

Therefore, we have obtained \textsc{Ramanujan}'s expansion into powers of $\frac{1}{m}$, and the coefficient of $\frac{1}{m^p}$ is \begin{align}
R_p=\left\{(-1)^{p-1}\frac{1}{2p8^p} + \sum_{s=0}^{p-1}\frac{D_{s}}{2^{s}}(-1)^{p-s}\binom{p-1}{p-s}\frac{1}{8^{p-s}}\right\}
\end{align}

But,\begin{align*}
\frac{D_{s}}{2^{s}}(-1)^{p-s}\binom{p-1}{p-s}\frac{1}{8^{p-s}}&=-\frac{\frac{B_{2s}\left(\frac{1}{2}\right)}{2s}}{2^{s}}(-1)^{p-s}\binom{p-1}{p-s}\frac{1}{8^{p-s}}\\
&=(-1)^{p-s-1}\frac{B_{2s}\left(\frac{1}{2}\right)}{2s2^{s}}\binom{p-1}{p-s}\frac{1}{8^{p-s}}\\
\end{align*}

\noindent and therefore\begin{align*}
R_p&=(-1)^{p-1}\frac{1}{2p8^p} + \sum_{s=0}^{p-1}\frac{D_{s}}{2^{s}}(-1)^{p-s}\binom{p-1}{p-s}\frac{1}{8^{p-s}}\\
&=(-1)^{p-1}\frac{1}{2p8^p} + \sum_{s=0}^{p-1}(-1)^{p-s-1}\frac{B_{2s}\left(\frac{1}{2}\right)}{2s2^{s}}\binom{p-1}{p-s}\frac{1}{8^{p-s}}\\
&=(-1)^{p-1}\left\{\frac{1}{2p8^p} + \sum_{s=1}^{p}(-1)^s\frac{B_{2s}\left(\frac{1}{2}\right)}{2s2^{s}}\binom{p-1}{p-s}\frac{1}{8^{p-s}}\right\}\\
&=(-1)^{p-1}\left\{\frac{1}{2p8^p} + \sum_{s=1}^{p}(-1)^s\frac{B_{2s}\left(\frac{1}{2}\right)}{2\cdot2^{s}}\cdot\frac{1}{p}\binom{p}{s}\frac{1}{8^{p-s}}\right\}\\
&=\frac{(-1)^{p-1}}{2p8^p}\left\{1 + \sum_{s=1}^{p}\binom{p}{s}(-4)^s B_{2s}\left(\frac{1}{2}\right)\right\}\\
\end{align*}

Thus. the formula for $H_n$ takes the form:
\begin{align}
H_n    &=\frac{1}{2}\ln(2m)+\gamma+\sum_{p=1}^{r}\frac{(-1)^{p-1}}{2p8^p}\left\{1 + \sum_{s=1}^{p}\binom{p}{s}(-4)^s B_{2s}\left(\frac{1}{2}\right)\right\}\cdot\frac{1} {m^{p}}\\
&\ \ \ +\epsilon_r+E_r+\theta_r\cdot\frac{D_{r+1}}{\left(n+\frac{1}{2}\right)^{2r+2}}   
      \end{align}We see\emph{ that (9) is the \textsc{Ramanujan} expansion with the general formula for the coefficient, $R_p$, of} $\frac{1}{m^p},$ as given in the statement of the theorem, while (10) is (an undeveloped form of) \emph{the error term. } 
      
      We will now estimate the error.

      To do so we will use the fact that the sum of a convergent alternating series, whose terms (taken with positive sign) decrease monotonically to zero, is equal to any partial sum \emph{plus a positive proper fraction of the first neglected term (with sign)}.
      
      Thus, $$\epsilon_r:=\sum_{l=r+1}^{\infty}(-1)^{l-1}\frac{1}{2l8^lm^l}=\alpha_{r}(-1)^{r}\frac{1}{2(r+1)8^{(r+1)}m^{r+1}}$$where $0<\alpha_{r}<1.$ 
      
      Moreover,\begin{align*}
      E_r&
      =\frac{D_{1}}{2^{1}}\sum_{k=r}^\infty(-1)^k\binom{k}{k}\frac{1}{8^k}\cdot\frac{1} {m^{1+k}}
+\frac{D_{2}}{2^{2}}\sum_{k=r-1}^\infty(-1)^k\binom{k+1}{k}\frac{1}{8^k}\cdot\frac{1} {m^{2+k}} 
+\cdots\\
&\ \ \ +\frac{D_{r}}{2^{r}}\sum_{k=1}^\infty(-1)^k\binom{k+r-1}{k}\frac{1}{8^k}\cdot\frac{1} {m^{r+k}} 
\\
&=\left\{\delta_1\frac{D_{1}}{2^{1}}(-1)^r\binom{r}{r}\frac{1}{8^r}+\delta_2\frac{D_{1}}{2^{2}}(-1)^{r-1}\binom{r}{r-1}\frac{1}{8^{r-1}}+\cdots+\delta_r\frac{D_{r}}{2^{r}}(-1)^1\binom{r}{1}\frac{1}{8^1}\right\}\frac{1}{m^{r+1}}\\      
&=\Delta_r\left\{\frac{D_{1}}{2^{1}}(-1)^r\binom{r}{r}\frac{1}{8^r}+\frac{D_{2}}{2^{2}}(-1)^{r-1}\binom{r}{r-1}\frac{1}{8^{r-1}}+\cdots+\frac{D_{r}}{2^{r}}(-1)^1\binom{r}{1}\frac{1}{8^1}\right\}
\end{align*}\noindent where $0<\delta_k<1$ for $k=1,\ 2, \cdots, \ r$ and $0<\Delta_r<1$.

Finally$$\theta_r\cdot\frac{D_{r+1}}{\left(n+\frac{1}{2}\right)^{2r+2}}=\theta_r\cdot\frac{D_{r+1}}{(2m)^{r+1}\left(1+\frac{1}{8m}\right)^{r+1}}=\delta_{r+1}\cdot\frac{D_{r+1}}{2^{r+1}}\cdot\frac{1}{m^{r+1}}$$\noindent where$0<\delta_{r+1}<1.$

 Thus, \emph{the total error is equal to} \begin{align*}&\epsilon_r+E_r+\theta_r\cdot\frac{D_{r+1}}{\left(n+\frac{1}{2}\right)^{2r+2}}\\
=&\Theta_r\cdot\left\{(-1)^{r}\frac{1}{2(r+1)8^{(r+1)}}+\sum_{q=1}^{r+1}\frac{D_{2q}}{2^{q}}(-1)^{r-q+1}\binom{r}{r-q+1}\frac{1}{8^{r-q+1}}\right\}\frac{1}{m^{r+1}}\\
=&\Theta_r\cdot R_{r+1} \end{align*}

\noindent by (8), where $0<\Theta_r<1,$ which is of the form as claimed in the theorem.  This completes the proof.

\end{proof}

\subsubsection*{Acknowledgment}
 Support from the Vicerrector\'{\i}a de Investigaci\'on of the 
University of Costa Rica is acknowledged.

\end{document}